\documentclass{article}
\usepackage{amsmath}
\usepackage{amssymb}

\usepackage{amscd}

\usepackage{amsthm}

\author{Diana Dziewa-Dawidczyk  %
 \and Zbigniew Pasternak-Winiarski\footnote{Faculty of Mathematics
and Information Science,  Warsaw University of Technology, Pl.
Politechniki 1, 00-661 Warsaw, POLAND}}
\title{Uniform structures  on differential spaces}

\begin{document}

\maketitle

\begin{abstract}
The uniform structure on a differential space defined by a family
of generators is considered. \bigskip


\end{abstract}

\section {Introduction}

This paper is the second of the series of publications concerning
integration of differential forms and densities on differential
spaces (the first one is \cite{ddpas}). It has the preliminary
character. We recall basic facts from the theory of differential
spaces and the theory of uniform structures. After that we
describe uniform structures defined on a differential space by
families of generators of its differential structure.

Section 2 of the paper contains basic definitions and the
description of preliminary facts concerning theory of differential
spaces. Foundations of theory of differential spaces can be find
in \cite{sik}. In Section 3 we give basic definitions and describe
the standard facts concerning theory of uniform spaces. We define
(in a standard manner) the uniform structure given on a
differential space  by a family of generators of its differential
structure. Section 4 contains basic facts concerning uniform
(uniformly continuous) maps. In Section 5 we recall the definition
of a complete uniform space  and the standard construction of
completion of a given uniform space. Here we introduce and
investigate the notion of the extension of a differential
structure.

 Without any other explanation we use the following symbols:
${\bf N}$-the set of natural numbers; ${\bf R}$-the set of reals.

\section {Differential spaces}

Let $M$ be a nonempty set and let $\mathcal{C}$ be a family of
real valued functions on $M$. Denote by $\tau_{\mathcal{C}}$ the
weakest topology on $M$ with respect to which all functions of
$\mathcal{C}$ are continuous.

A base of the topology $\tau_{\mathcal C}$ consists of sets:\\
\begin{displaymath}
(\alpha_1,\ldots,\alpha_n)^{-1}(P)=\bigcap_{i=1}^{n} \{p:
a_i<\alpha_i(p)<b_i\},
\end{displaymath}
where $n\in{\bf N}$, $a_1,\ldots,a_n,b_1,\ldots,b_n\in{\bf R}$,
$a_i<b_i$, $\alpha_1, \ldots,\alpha_n\in{\mathcal C}$,
$P=\{(x_1,\ldots,x_n)\in{\bf R}^{n};a_i<x_i<b_i,i=1,
\ldots,n\}$.\\

{\sc Definition 2.1} A function $f: M\rightarrow \bf{R}$ is called
\emph{a local ${\mathcal C}$-function on} $M$ if for every $p\in M$
there is a neighbourhood $V$ of $p$ and $\alpha\in {\mathcal C}$
such that $f_{|V}=\alpha_{|V}$. The set of all local
${\mathcal C}$-functions on $M$ is denoted by ${\mathcal C}_M$.\\

Note that any function $f\in{\mathcal C}_M$ is continuous with
respect to the topology $\tau_{\mathcal C}$. In fact, if
$\{V_i\}_{i\in I}$ is such an open (with respect to
$\tau_{\mathcal C}$) covering of $M$ that for any $i\in I$ there
exists $\alpha_i\in{\mathcal C}$ satisfying
$f_{|V_i}=\alpha_{i|V_i}$ and $U$ is an open subset of ${\bf R}$
then
$$f^{-1}(U)=\bigcup\limits_{i\in I}(\alpha_{i|V_i})^{-1}(U).$$
Since $(\alpha_{i|V_i})^{-1}(U)$ is open in $V_i$ and $V_i\in
\tau_{\mathcal C}$ we obtain $(\alpha_{i|V_i})^{-1}(U)\in
\tau_{\mathcal C}$ for any $i\in I$. Hence $f^{-1}(U)\in
\tau_{\mathcal C}$. Bearing in mind that $U$ is an arbitrary open
set in ${\bf R}$ we obtain that $f$ is continuous with respect to
$\tau_{\mathcal C}$.

 We have ${\mathcal C}\subset{\mathcal C}_M$ which implies
 $\tau_{\mathcal C}\subset\tau_{{\mathcal C}_M}$. On the other
 hand any element of ${\mathcal C}_M$ is a function continuous
 with respect to $\tau_{\mathcal C}$. Then $\tau_{{\mathcal C}_M}\subset
 \tau_{\mathcal C}$ and consequently $\tau_{{\mathcal C}_M}=
 \tau_{\mathcal C}$.\\

{\sc Definition 2.2} A function $f: M\rightarrow \bf{R}$ is called
\emph{${\mathcal C}$-smooth function on} $M$ if there exist
$n\in{\bf N}$, $\omega \in C^{\infty} ({\bf R}^{n})$ and
$\alpha_1,\ldots,\alpha_n\in{\mathcal C}$ such that
$$f=\omega \circ(\alpha_1,\ldots,\alpha_n).$$
The set of all ${\mathcal C}$-smooth functions on $M$ is denoted
by ${\rm sc}{\mathcal C}$.\\

We have ${\mathcal C}\subset sc{\mathcal C}$ which implies
$\tau_{\mathcal C}\subset\tau_{sc{\mathcal C}}$. On the other hand
any superposition $\omega \circ(\alpha_1,\ldots,\alpha_n)$ is
continuous with respect to $\tau_{\mathcal C}$ which gives
$\tau_{sc{\mathcal C}}\subset\tau_{\mathcal C}$. Consequently
$\tau_{sc{\mathcal C}}=\tau_{\mathcal C}$.\\

{\sc Definition 2.3} A set $\mathcal{C}$ of real functions on $M$
is said to be a \emph{(Sikorski's) differential structure} if: (i)
$\mathcal{C}$ is \emph{closed with respect to localization} i.e.
${\mathcal C} = {\mathcal C}_M$; (ii) ${\mathcal C}$ is {\sl
closed with respect to superposition with smooth functions} i.e.
${\mathcal C}={\rm sc}{\mathcal C}$.\\

In this case a pair $(M,{\mathcal C})$ is said to be
\emph{a (Sikorski's) differential space}.\\

{\sc Proposition 2.1.} {\it The intersection of any family of
differential structures defined on a set $M\ne\emptyset$ is a
differential structure on $M$.

\smallskip
Proof}. Let $\{{\mathcal C}_i\}_{i \in I}$ be a family of
differential structures defined on a set $M$ and let ${\mathcal
C}:=\bigcap_{i \in I}{\mathcal C}_i$. Then ${\mathcal C}$ is
nonempty family of real-valued functions on $M$ (it contains all
constant functions). If $n\in{\bf N}$, $\omega \in C^{\infty}
({\bf R}^{n})$ and $\alpha_1,\ldots,\alpha_n\in{\mathcal C}$ then
for any $i\in I$ $\alpha_1,\ldots,\alpha_n\in{\mathcal C}_i$ and
consequently $\omega \circ(\alpha_1,\ldots,\alpha_n)\in{\mathcal
C}_i$. Hence $\omega\circ(\alpha_1,\ldots,\alpha_n)\in{\mathcal
C}$ which means that $sc{\mathcal C}={\mathcal C}$.

Since ${\mathcal C}\subset{\mathcal C}_i$ for any $i\in I$ we have
$\tau_{\mathcal C}\subset\tau_{{\mathcal C}_i}$. It means that any
subset of $M$ open with respect to $\tau_{\mathcal C}$ is open
with respect to $\tau_{{\mathcal C}_i}$, for $i\in I$.

Let $\beta\in{\mathcal C}_M$. Choose for any $p\in M$ a set
$U_p\in\tau_{\mathcal C}$ and a function $\alpha_p\in{\mathcal C}$
such that $p\in U_p$ and $\beta_{|U_p}=\alpha_{p|U_p}$. Since
$\alpha_p\in{\mathcal C}_i$ and $U_p\in\tau_{{\mathcal C}_i}$ we
obtain $\beta\in({\mathcal C}_i)_M={\mathcal C}_i$, for any $i\in
I$. Then $\beta\in{\mathcal C}$ and consequently ${\mathcal
C}_M={\mathcal C}$.

Equalities ${\mathcal C}_M={\mathcal C}=sc{\mathcal C}$ means that
${\mathcal C}$ is a differential structure on $M$.\hfill $\Box$\\

Let ${\mathcal F}$ be a set of real functions on $M$. Then, by
Proposition 2.1, the intersection ${\mathcal C}$ of all
differential structures on $M$ containing ${\mathcal F}$ is a
differential structure on $M$. It is the smallest differential
structure on $M$ containing ${\mathcal F}$. One can easy prove
that ${\mathcal C}=(sc{\mathcal F})_M$ (see \cite{wal}). This
structure is called \emph{the differential structure generated by}
${\mathcal F}$. Functions of ${\mathcal F}$ are called
\emph{generators} of the differential structure ${\mathcal C}$. We
have also $\tau_{(sc{\mathcal F})_M}=\tau_{sc{\mathcal
F}}=\tau_{\mathcal F}$
(see remarks after Definitions 2.1 and 2.2).\\

\smallskip
Let $(M,{\mathcal C})$ and $(N, {\mathcal D})$ be differential spaces.
A map $F: M \rightarrow N$ is said to be \emph{smooth} if for any
$\beta\in{\mathcal D}$ the superposition $\beta\circ F\in{\mathcal C}$.
We will denote the fact that $f$ is smooth writing
$$F:(M,{\mathcal C})\to(N,{\mathcal D}).$$
If $F:(M,{\mathcal C})\to(N,{\mathcal D})$ is a bijection and
$F^{-1}:(N,{\mathcal D})\to(M,{\mathcal C})$ then $F$ is called
\emph{a diffeomorphism}

It is easy to show that if $\mathcal{F}$ is a family of generators
of the structure ${\mathcal D}$ on the set $N$ then
$F:(M,{\mathcal C})\to(N,{\mathcal D})$ iff for any
$f\in\mathcal{F}$ the superposition $f\circ F\in\mathcal{C}$.

\smallskip
If $A$ is a nonempty subset of $M$ and ${\mathcal C}$ is a differential
structure on $M$ then ${\mathcal C}_A$ denotes the differential structure
on $A$ generated by the family of restrictions
$\{\alpha_{|A}:\alpha\in{\mathcal C}\}$. The differential space
$(A,{\mathcal C}_A)$ is called \emph{a differential subspace} of
$(M,{\mathcal C})$. One can easy prove the following\\

{\sc Proposition 2.2.} {\it Let $(M,{\mathcal C})$ and $(N, {\mathcal D})$
be differential spaces and let $F:M\to N$. Then $F:(M,{\mathcal C})\to
(N,{\mathcal D})$ iff} $F:(M,{\mathcal C})\to(F(M),F(M)_{\mathcal D}).$

\medskip
If the map $F:(M,{\mathcal C})\to(F(M),F(M)_{\mathcal D})$ is a
diffeomorphism then we say that $F:M\to N$ is \emph{a
diffeomorphism onto its range} (in $(N,{\mathcal D}))$. In
particular the natural embedding $A\ni x\mapsto i(x):= x\in M$ is
a diffeomorphism of $(A,{\mathcal C}_A)$ onto its range in
$(M,{\mathcal C})$.

\bigskip
If $\{(M_i,{\mathcal C}_i)\}_{i\in I}$ is an arbitrary family of differential
spaces then we consider the Cartesian product $\prod\limits_{i\in I}M_i$ as
a differential space with the differential structure
$\bigotimes\limits_{i\in I}{\mathcal C}_i$ generated by the family of
functions ${\mathcal F}:=\{\alpha_i\circ pr_i:i\in I,\alpha_i\in{\mathcal C}_i\}$,
where $\prod\limits_{i\in I}M_i\ni(x_i)\mapsto pr_j((x_i))=:x_j\in M_j$ for
any $j\in I$. The topology $\tau_{\bigotimes\limits_{i\in I}{\mathcal C}_i}$
coincides with the standard product topology on $\prod\limits_{i\in I}M_i$.

\bigskip
\emph{A generator embedding} of the differential space
$(M,{\mathcal C})$ into the Cartesian space is a mapping
$\phi_{{\mathcal F}}: (M,{\mathcal C}) \rightarrow
(\textbf{R}^{{\mathcal F}},C^{\infty} (\textbf{R}^{{\mathcal
F}}))$ given by the formula
$$\phi_{{\mathcal F}}(p)=(\alpha(p))_{\alpha \in {\mathcal F}}$$
(for example if ${\mathcal F} = \{\alpha_1,\alpha_2,\alpha_3\}$
then $\phi_{{\mathcal
F}}(p)=(\alpha_1(p),\alpha_2(p),\alpha_3(p))\in \textbf{R}^3
\cong\textbf{R}^{\mathcal F}$).

\smallskip
{\sc Proposition 2.3.} {\it Let ${\mathcal F}$ be a family of
generators of the differential structure ${\mathcal C}$ on the set
$M$ If ${\mathcal F}$ separates points of $M$ then the generator
embedding defined by ${\mathcal F}$ is a diffeomorphism onto its
image. On that image we consider a differential structure of a
subspace of $(\textbf{R}^{\mathcal F},C^{\infty}
(\textbf{R}^{{\mathcal F}}))$.

\smallskip
Proof}. Since ${\mathcal F}$ separates points of $M$ it follows
from the definition of differential embedding $\phi_{{\mathcal
F}}$ that it is an one-to-one mapping onto its image in
$\textbf{R}^{{\mathcal F}}$. Moreover for any $f\in{\mathcal F}$
we have $pr_f\circ\phi_{{\mathcal F}}=f\in{\mathcal C}$. Since the
differential structure $C^{\infty} (\textbf{R}^{{\mathcal F}})$ is
generated by the family $\{pr_g\}_{g\in{\mathcal C}}$ we obtain
that the map $\phi_{{\mathcal F}}$ is smooth with respect to
$C^{\infty} (\textbf{R}^{{\mathcal F}})$. It remains to show that
the map $\phi_{{\mathcal F}}^{-1}$ is smooth.

For any $f\in{\mathcal F}$ we have
$$f\circ\phi_{{\mathcal F}}^{-1}=pr_{f|{{\mathcal F}}(M)}.$$
It means that $f\circ\phi_{{\mathcal F}}^{-1}\in C^{\infty}
(\textbf{R}^{{\mathcal F}})_{{{\mathcal F}}(M)}$. Since the
differential structure $\mathcal{C}$ is generated by ${\mathcal
F}$ we obtain $\phi_{{\mathcal F}}^{-1}$ is smooth. $\Box$

\section {Uniform structures}

Let $X$ be a nonempty set.\\

{\sc Definition 3.1.} A set $\Delta = \{(x,x):x\in X\}$ is said to
be {\sl the diagonal} of the product $X\times X$. A set $V\subset
X\times X$ is called {\sl a neighbourhood of the
diagonal} if $\Delta\subset V$ and $V=-V$, where $-V=\{(x,y): (y,x)\in V\}$.\\
A family of all neighborhoods of the diagonal is denoted by ${\mathcal D}_{X}$.\\

{\sc Definition 3.2} If for $x,y\in X$ and $V\in {\mathcal D}_{X}$
we have $(x,y)\in V$, then we say that {\sl $x$ and $y$ are
distant less then $V$} $(|x-y|<V)$. We say that {\sl the diameter
of a set $A\subset X$ is less then $V$} $(\delta(A)<V)$ if for all
$x,y\in A$ we have $|x-y|<V$. {\sl A ball with the center at
$x_0\in X$ and the radius} $V$ is a set $K(x_0,V)=\{x\in
X:|x_0-x|<V\}.$ The set
$$2V:=\{(x,y)\in X\times X:\exists z\in X\ [(x,z)\in V\wedge(x,z)\in V]\}.$$

{\sc Definition 3.3}  {\sl A uniform structure} $\mathcal{U}$ on $X$ is a subfamily of
${\mathcal D}_{X}$ satisfying the following conditions:\\

1) $(V\in{\mathcal U} \wedge V\subset W\in{\mathcal D}_{X}\Rightarrow (W\in{\mathcal U})$;\\

2) $(V_{1},V_{2}\in{\mathcal U})\Rightarrow(V_{1}\cap V_{2}\in{\mathcal U})$;\\

3) $\forall V\in{\mathcal U}\ \exists W\in{\mathcal U}\ [2W\subset V]$;\\

4) $\bigcap{\mathcal U}=\Delta$.\\

If ${\mathcal U}$ is a uniform structure on $X$ then the pair
$(X,\mathcal{U})$ is called
{\sl a uniform space}.\\

{\sc Definition 3.4} {\sl A base} of a uniform structure ${\mathcal U}$ in $X$ is a family
${\mathcal B}\subset{\mathcal U}$ such, that for all $V\in{\mathcal U}$ there exists
$W\in{\mathcal B}$ satisfying $W\subset V$.\\

Each base ${\mathcal B}$ has following properties:\\

B1) $(V_{1},V_{2}\in{\mathcal B})\Rightarrow (\exists V\in{\mathcal B}\ [V\subset
V_{1}\cap V_{2}])$;\\

B2) $\forall V\in{\mathcal B}\ \exists W\in{\mathcal B}\ [2W\subset V]$;\\

B3) $\bigcap{\mathcal B}=\Delta$.\\

On the other hand it can be easy proved that if a family
${\mathcal B}$ of neighbourhoods of the diagonal of a set $X$
fulfiels conditions (B1)-(B3) then there exists exactly one
uniform structure ${\mathcal U}$ on $X$ such that ${\mathcal B}$
is a base of ${\mathcal U}$.

\smallskip
Every neighbourhood $V\in {\mathcal D}_{X}$ of the diagonal
defines the covering ${\mathcal P}(V)= \{K(x,V)\}_{x\in X}$ of the
set X. If ${\mathcal U}$ is a uniform structure in $X$ then every
covering ${\mathcal O}$ of $X$ for which there exists
$V\in{\mathcal U}$ such that ${\mathcal P}(V)$ is a refinement of
${\mathcal O}$ is said to be {\sl a uniform
covering} (with respect to ${\mathcal U}$).\\

Each uniform structure on X defines a topology on $X$. In other words each uniform space
$(X,{\mathcal U})$ defines a topological space $(X,\Theta)$.\\

{\sc Theorem 3.1} {\it If ${\mathcal U}$ is a uniform structure on $X$, then a family
\linebreak
$\Theta=\{G\subset X:\forall x\in G\quad \exists V\in{\mathcal U}\quad [K(x,V)\subset G]\}$
is a topology in X and $(X,\Theta)$ is $T_{1}$-space. A topology $\Theta$ is said to be
{\sl a topology given in $X$ by uniform structure} ${\mathcal U}$ and is denoted by
$\tau_{\mathcal U}$}.\\

\smallskip
For the proof see \cite{Bour} or \cite{eng}.

\bigskip
It can be proved that a topology $\tau$ on a topological space $X$ is given by some uniform
structure on $X$ if and only if $X$ is a Tichonov space (see \cite{eng}).\\

Let $\varrho$ be a pseudometric on a uniform space $(X,{\mathcal U})$. If for every
$\varepsilon > 0$ there is $V\in{\mathcal U}$ such that if $|x-y|<V$ then $\varrho (x,y) <
\varepsilon$, then $\varrho$ is called {\sl a uniform pseudometric (with respect to
$\mathcal{U}$)}.\\

We can defined a uniform structure on three different ways: (i) if we give a
base; (ii) if we give a family of uniform coverings or (iii) if we give a family of
pseudometrics (see \cite{eng}).\\

Let $(M,{\mathcal C})$ be a differential space such that $\mathcal{C}={\rm sc}{\mathcal F}_M$
and $(M,\tau_{\mathcal C})$ is a Hausdorff space (the last is true iff the family
${\mathcal C}$ separates points in $X$ iff the family ${\mathcal F}$ separates points in $X$).
On the set $M$ the family ${\mathcal F}$ defines the uniform structure
${\mathcal U}_{\mathcal F}$ such that the base ${\mathcal B}$ of ${\mathcal U}_{\mathcal F}$
is given as follows:\\
$${\mathcal B}=\{V(f_{1},\ldots,f_{k},\varepsilon)\subset M\times M;
k\in{\bf N};f_{1},\ldots,f_{k}\in{\mathcal F}, \varepsilon>0\},$$
where
$$V(f_{1},\ldots,f_{k},\varepsilon)=\{(x,y)\in M\times M:\forall
1\leq i\leq k\quad |f_{i}(x)-f_{i}(y)|<\varepsilon\}.$$

{\sc Proposition 3.1} {\it The family ${\mathcal B}$ satisfies on $M$ conditions B1 - B3.

\smallskip
Proof.} (B1) Let: $V_1=V(f_1,\ldots,f_k,\varepsilon_1)\in
{\mathcal B},\ V_2=V(g_{1},\ldots,g_{n},\varepsilon_2)\in
{\mathcal B}$ and $\varepsilon={\mathrm min}(\varepsilon_1, \varepsilon_2)$. Then
$$V:=V(f_{1},\ldots,f_{k},g_{1},\ldots,g_{n},\varepsilon)=$$
$$\{(x,y)\in M\times M: \forall 1\leq i\leq k\ [|f_{i}(x)-f_{i}(y)|<\varepsilon]\wedge
\forall 1\leq j\leq n\ [|g_{i}(x)-g_{i}(y)|<\varepsilon]\}\in{\mathcal B}$$
and $V\subset V_{1}\cap V_{2}$.

(B2) Let $V=V(f_{1},\ldots,f_{k},\varepsilon)\in \mathcal{B}$. Then $W:=
V(f_{1},\ldots,f_{k},\frac{\varepsilon}{2})\in{\mathcal B}$ and
$$2W=$$
$$\{(x,y)\in M\times M:\exists z\in M\ \forall 1\leq i\leq k\ [|f_{i}(x)-f_{i}(z)|<
\frac{\varepsilon}{2}\ \wedge\ |f_{i}(z)-f_{i}(y)|<\frac{\varepsilon}{2}]\}$$
$$\subset\{(x,y)\in M\times M:\forall 1\leq i\leq k\ [|f_{i}(x)-f_{i}(y)|<\varepsilon]\}=V.$$

(B3) Since for any $V\in{\mathcal B}$ there is $\Delta\subset V$ we have
$$\Delta\subset\bigcap{\mathcal B}.$$
On the other hand
$$\bigcap{\mathcal B}\ \ \subset\bigcap\limits_{f\in{\mathcal F},\varepsilon>0}V(f,\varepsilon)=$$
$$\{(x,y)\in M \times M:\forall f\in{\mathcal F}\ \ \forall \varepsilon>0\quad
[|f(x)-f(y)|<\varepsilon]\}=$$
$$\{(x,y)\in M \times M:\forall f\in{\mathcal F}\quad [|f(x)=f(y)]\}=$$
$$\{(x,x)\in M\times M\} = \Delta.\ \ \Box$$

\smallskip
The uniform space $(M,{\mathcal U}_{\mathcal F})$ is said to be {\sl the uniform
space given by the family of function} ${\mathcal F}$.

\bigskip
If we have two different families ${\mathcal F}_1$, and ${\mathcal F}_2$ of generators of
differential space $(M,{\mathcal C})$ then the uniform structures ${\mathcal U}_{{\mathcal F}_1}$
and ${\mathcal U}_{{\mathcal F}_1}$ can be different too.

\smallskip
{\sc Example 3.1} Let $M={\bf R}$, ${\mathcal C}=C^\infty({\bf R})$, ${\mathcal F}_1=
\{id_{\bf R}\}$ and ${\mathcal F}_2=\{id_{\bf R},f\}$, where
$$id_{\bf R}(x)=x,\quad{\mathrm and}\quad f(x)=x^2,\qquad x\in{\bf R}.$$
Then does not exists $\varepsilon >0$ such that $V(id_{\bf
R},\varepsilon)\subset V(f,1)$. Hence $V(f,1)\notin {\mathcal
U}_{{\mathcal F}_1}$ and ${\mathcal U}_{{\mathcal F}_1}\neq
{\mathcal U}_{{\mathcal F}_2}$. $\ \ \Box$\\

\section {Uniform continuous mapping}

Let $(X,{\mathcal U})$, $(Y,{\mathcal V}), (X,{\mathcal \overline{U}})$,
$(Y,{\mathcal \overline{V}})$ be uniform spaces.\\

\smallskip
{\sc Definition 4.1}  A mapping $f:X\to Y$ is said to be {\sl uniform} with respect to
uniform structures ${\mathcal U}$ and ${\mathcal V}$ if
$$\forall V\in{\mathcal V}\ \exists U\in{\mathcal U}\ \forall x,x'\in X\
[|x-x'|<U\Rightarrow |f(x)-f(x')|]<V.$$ In other words for every
$V\in{\mathcal V}$ there is $U\in{\mathcal U}$ such that $U\subset
(f\times f)^{-1}(V)$. We denote it by
$$f:(X,{\mathcal U})\to(Y,{\mathcal V}).$$

It is easy to prove that:\\
(i) any uniform mapping $f:(X,{\mathcal U})\to(Y,{\mathcal V})$ is
continuous with respect to topologies $\tau_{\mathcal U}$ and $\tau_{\mathcal V}$;\\
(ii) a superposition of uniform mappings is a uniform mapping.\\

\bigskip
We can give criteria of the uniformity:\\

\smallskip
{\sc Theorem 4.1} {\it Let $f:X\to Y$ and let ${\mathcal U}$ and ${\mathcal V}$ be
uniform structures on $X$ and $Y$ respectively. Then the following conditions are
equivalent:\\

(a) $f:(X,{\mathcal U})\to(Y,{\mathcal V})$.\\

(b) If ${\mathcal B}$ and ${\mathcal D}$ are bases of ${\mathcal U}$ and ${\mathcal V}$
respectively then for each $V\in{\mathcal D}$ there exists $U\in{\mathcal B}$ such that
$U\subset (f\times f)^{-1}(V)$.\\

(c) For every covering $\mathcal{A}$ of $Y$ uniform with respect to ${\mathcal V}$,
a covering $\{f^{-1}(A)\}_{A\in{\mathcal A}}$ of $X$ is uniform with respect to
${\mathcal U}$.\\

(d) For every pseudometric $\varrho$ on $Y$ uniform with respect
to ${\mathcal V}$, a pseudometric $\sigma$ on $X$ given by the
formula
$$\sigma (x,y)=\varrho(f(x),f(y))\qquad x,y\in X$$
is uniform with respect to the uniform structure} ${\mathcal U}$.

\smallskip
For the proof see \cite{eng}.

\bigskip
A mapping $f$, that is a uniform with respect to uniform
structures ${\mathcal U}$ and ${\mathcal V}$ could be not uniform
with respect to uniform structures $\overline{{\mathcal U}}$ and
$\overline{{\mathcal V}}$.\\

Example 4.1. Let $M=\textbf{R}$, $C=C^{\infty}(\bf{R})$,
${\mathcal F}_{1} = \{id_{\bf{R}}\}, {\mathcal F}_{2} =
\{id_{\bf{R}},f\}$, where $f(x)=x^{2}$, $x \in \bf{R}$.

Here $f$ is the uniform mapping with respect to ${\mathcal
F}_{2}$, but it is not uniform with respect to ${\mathcal F}_{1}$.
In fact, the set $V = \{(x,y): |f(x) - f(y)| = x^2 - y^2 <
\varepsilon \}$ is an element of ${\mathcal D}$ and does not
exists $U \in {\mathcal B}$ such that ${\mathcal U} \subset (f
\times f)^{-1}(V)$ (see Example 3.1).\\

\bigskip
{\sc Definition 4.2} A bijective mapping $f:(X,{\mathcal U})\to(Y,{\mathcal V})$ is
{\sl a uniform homeomorphism} if $f^{-1}$ is a uniform mapping. Then we say that
$(X,{\mathcal U})$ and $(Y,{\mathcal V})$ are {\sl uniformly homeomorphic}.

\smallskip
By (i) it is obvious that if $f:(X,{\mathcal U})\to(Y,{\mathcal
V})$ is a uniform homeomorphism then $f$ is a homeomorphism of the
topological spaces $(X,\tau_{\mathcal U})$ and $(Y,\tau_{\mathcal
V})$.\\

\section {Complete uniform spaces and extensions of differential
structure}

{\sc Definition 5.1} Let $X$ be a non empty set, $x\in X$ and
$V\in{\mathcal D}_{X}$ (see Definition 3.1). A set $U\subset X$ is
said to be {\sl small of rank} $V$ if
$\exists x\in V\ [U\subset K(x,V)]$ (see Definition 3.2).\\

{\sc Definition 5.2} A nonempty family ${\mathcal F}$ of subsets of a set
$X$ is said to be a filter on $X$ if:\\

(F1) $(F\in{\mathcal F}\ \wedge\ F\subset U\subset X)\ \Rightarrow\ (U\in{\mathcal F})$;\\

(F2) $(F_1,F_2 \in{\mathcal F})\Rightarrow(F_1 \cap F_2 \in{\mathcal F})$;\\

(F3) $\emptyset\notin{\mathcal F}$.\\

{\sc Definition 5.3} A filtering base on $X$ is a nonempty family
${\mathcal B}$ of subsets of $X$ such that\\

(FB1) $\forall A_1,A_2 \in{\mathcal B}\ \exists A_3\in{\mathcal B}\
[A_3\subset A_1\cap A_2]$;\\

(FB2) $\emptyset\notin{\mathcal B}$.\\

\smallskip
If ${\mathcal B}$ is a filtering base on $X$ then
$${\mathcal F}=\{F\subset X:\exists A\in{\mathcal B}\ [A\subset F]\}$$ is a filter
on $X$. It is called {\sl the filter defined by} ${\mathcal B}$.\\

\smallskip
{\sc Definition 5.4} Let $X$ be a topological space. We say that a filter $\mathcal{F}$
on $X$ {\sl is convergent to} $x\in X$ (${\mathcal F} \rightarrow x$) if  for any
neighbourhood $U$ of $x$ there exists $F\in{\mathcal F}$ such that $F \subset U$.\\

\smallskip
{\sc Definition 5.5} Let $(X,{\mathcal U})$ be a uniform space. A filter ${\mathcal F}$
on $X$ is {\sl a Cauchy filter} if
$$\forall V\in{\mathcal U}\ \ \exists F\in{\mathcal F}\ \ [F\times F\subset V].$$

\smallskip
{\sc Definition 5.6}  A uniform space $(X,{\mathcal U})$ is said to be {\sl complete} if
each Cauchy filter on $X$ is convergent in $\tau_{\mathcal U}$.

\bigskip
Let $(X,{\mathcal U})$ be a uniform space, $M\subset X$ and $M
\neq \O$. Denote
$${\mathcal U}_M:=\{V\cap M:V\in{\mathcal U}\}.$$
Then it is easy to show that ${\mathcal U}_M$ is a uniform structure on $M$. We call
$(M,{\mathcal U}_M)$ {\sl a uniform subspace of the uniform space} $(X,{\mathcal U})$.

\bigskip
{\sc Theorem 5.1} {\it If $(X,{\mathcal U})$ is a complete uniform
space and $M$ is a closed subset of the topological space
$(X,\tau_{\mathcal U})$ then a uniform space $(M,{\mathcal U}_M)$
is complete. Conversely, If $(M,{\mathcal U}_M)$ is a complete
uniform subspace of some (not necessarily complete) uniform space
$(X,{\mathcal U})$ then M is closed in $X$ with respect to}
$\tau_{\mathcal U}$.

\smallskip
For the proof see \cite{Bour}, \cite{eng} or \cite{pas}.

\bigskip
Any uniform space can be treated as a uniform subspace of some
complete uniform space. We have the following

\smallskip
{\sc Theorem 5.2} {\it  For each uniform space} $(X,{\mathcal U})$:\\
(i) {\it there exists a complete uniform space
$(\widetilde{X},\widetilde{\mathcal U})$ and a set
$A\subset\widetilde{X}$ dense in $\widetilde{X}$ (with respect to
the topology $\tau_{\widetilde{\mathcal U}}$) such that
$(X,\mathcal{U})$ is uniformly homeomorphic to}
$(A,\widetilde{\mathcal U}_A)$;\\
(ii) {\it if the complete uniform spaces $(\widetilde{X}_1,\widetilde{\mathcal U}_1)$
and $(\widetilde{X}_2,\widetilde{\mathcal U}_2)$ satisfies condition of the point} (i)
{\it then they are uniformly homeomorphic}.

\smallskip
For the details of the proof see \cite{Bour} or \cite{pas}. Here
we only want to describe the construction of
$(\widetilde{X},\widetilde{\mathcal U})$.

Let $\widetilde{X}$ be the set of all minimal (with respect to the
order defined by inclusion) Cauchy filters in $X$. For every
symmetric set $V\in{\mathcal U}$ we denote by $\widetilde{V}$ the
set of all pairs $({\mathcal F}_1,{\mathcal F}_2)$ of minimal
Cauchy's filters, which have a common element being a small set in
rank V. We define a family $\widetilde{\mathcal U}$ of subsets of
set $\widetilde{X}\times \widetilde{X}$ as the smallest uniform
structure on $X$ containing all sets from the family
$\{\widetilde{V}:V\in{\mathcal U}\}$.\\

\bigskip
Let us consider two different uniform structures at the same
differential space $(\textbf{R},{\mathcal C}^{\infty})$:
${\mathcal U}_{\mathcal{F}}$ and  ${\mathcal U}_{\mathcal{G}}$,
where ${\mathcal F}=\{id_{\bf{R}}\}$, ${\mathcal G} =\{arctg x\}$.
Then $(\textbf{R},{\mathcal U}_{\mathcal{F}}$ is the complete
space ($(\widetilde{\bf R} = {\bf R}$) whereas $(\textbf{R},
{\mathcal U}_{\mathcal{G}}$. In this case we can identify
$\widetilde{\bf R}$ with the interval $[-\frac{\pi}{2} ;
\frac{\pi}{2}]$.

\bigskip
Let $N$ be a set, $M \subseteq N$, $M \neq \O$, ${\mathcal C}$ be
a differential structure on $M$.

\smallskip
{\sc Definition 5.7.} The differential structure ${\mathcal D}$ on
$N$ is an \emph{extension} of the differential structure
${\mathcal C}$ from the set $M$ to the set $N$ if ${\mathcal C} =
{\mathcal D}_{M}$ (if we get the structure ${\mathcal C}$ by
localization of the structure ${\mathcal D}$ to $M$).

\bigskip
For the sets $N,M$ and the differential structure ${\mathcal C}$
on $M$ we can construct many different extensions of the structure
$M$ to $N$.

\smallskip
{\sc Example 5.1.} If for each function $f \in {\mathcal C}$ we
assign the function $f_{0} \in {\mathbf R}^{N}$ such that $f_{0|M}
= f$ and $f_{0|N \backslash M} \equiv 0$. Then the differential
structure generated on $N$ by the family of functions
$\{f_{0}\}_{f \in {\mathcal C}}$ is the extension of ${\mathcal
C}$ from $M$ to $N$. Similarly, if for each function $f \in
{\mathcal C}$ we assign the family of the functions ${\mathcal
F}_{f} := \{g \in {\mathbf R}^{N} : g_{|M} = f\}$, then the
differential structure on $N$ generated the family of the
functions ${\mathcal F} := \bigcup_{f \in {\mathcal C}} {\mathcal
F}_{f}$ is the extension of ${\mathcal C}$ from $M$ to $N$. If the
set $N \backslash M$ contains at least two elements, then the
differential structures generated by the families $\{f_{0}\}_{f
\in {\mathcal C}}$ and ${\mathcal F}$ are different.

\bigskip
{\sc Definition 5.8.} If $\tau$ is a topology on the set $N$, then
the extension ${\mathcal D}$ of the differential structure
${\mathcal C}$ from $M$ to $N$ is \emph{continuous with respect to
$\tau$} if each function $f \in {\mathcal D}$ is continuous in the
topology $\tau$ $(\tau_{\mathcal D}\subset\tau)$.

\bigskip
If on the set $N$ there exists continuous (with respect to $\tau$)
extension of the differential structure ${\mathcal C}$ from the
set $M \subset N$, then the structure ${\mathcal C}$ is said to be
 \emph{extendable from the set $M$ to the topological space} $(N, \tau)$.

\bigskip
{\sc Example 5.2.} The differential structure
$C^{\infty}(\textbf{R})_{\textbf{Q}}$ is extendable from the set
of rationales to the set of reals. The continuous extensions are
e.g. $C^{\infty}(\textbf{R})$ and the structure ${\mathcal D}$
generated on \textbf{R} by the family of functions
$C^{\infty}(\textbf{R}) \cup\{f\}$, where $f : \textbf{R}
\rightarrow \textbf{R}$, $f(x) := |x - \sqrt{2}|$, $x \in
\textbf{R}$.\\

It is not difficult to show that if ${\mathcal F}$ is a family of
generators of a differential structure ${\mathcal C}$ on a set $M$
then the completion $\widetilde{M}$ of $M$ with respect to the
uniform structure ${\mathcal U}_{\mathcal F}$ can be identify with
the closure of the range $\phi_{{\mathcal F}}(M)$ of the generator
embedding $\phi_{{\mathcal F}}$ in the Cartesian product ${\bf
R}^{\mathcal F}$. In this case the differential structure
$C^{\infty} (\textbf{R}^{\mathcal F})_{\phi_{{\mathcal F}}(M)}$ is
a natural continuous extension of ${\mathcal C}$ from $M$ to
$\widetilde{M}$.

\vspace{1cm}

\footnotesize
\begin{center}

\end{center}


\begin{thebibliography}{99}

\bibitem{Bour} N. Bourbaki, {\sl General Topology}, Springer-Verlag,
Berlin Heidelberg New York London Paris Tokyo 1989.

\bibitem{eng} R. Engelking, {\sl Topologia og\'olna II}, PWN, Warsaw 1989

\bibitem{pas} Z. Pasternak-Winiarski, {\sl Grupowe struktury
r\'o\.zniczkowe i ich podstawowe w\l asno\'sci (doctor thesis)},
Warsaw University of Technology, Warsaw 1981

\bibitem{wal} W. Waliszewski, {\sl Regular and coregular mappings of
differential space}, Annales Polonici Mathematici XXX, 1975

\bibitem{sik} R. Sikorski, {\sl Wst\c{e}p do geometrii r\'o\.zniczkowej},
PWN, Warsaw 1972

\bibitem{ddpas} D. Dziewa-Dawidczyk, Z. Pasternak-Winiarski, {\sl Differential structures on natural
bundles connected with a differential space}, "Singularities and
Symplectic Geometry VII" Singularity Theory Seminar (2009), S.
Janeczko (ed), Faculty of Mathematics and Information Science,
Warsaw University of Technology.
\end{thebibliography}
\end{document}